\documentclass[12pt]{amsart}
\usepackage{mathpazo}
\usepackage{charter}
\usepackage{amsfonts, amsmath, amssymb, lineno}
\usepackage{amsmath}
\usepackage{amssymb}
\usepackage{amsfonts}
\usepackage{graphicx}

\theoremstyle{plain}

\newtheorem{lemma}{Lemma}

\newtheorem{remark}{Remark}

\newtheorem{theorem}{Theorem}
\numberwithin{equation}{section}

\begin{document}
\title[The metric dimension of regular bipartite graphs]{The metric
dimension of regular bipartite graphs}
\author{S.W. Saputro}
\address[S.W. Saputro]{Combinatorial Mathematics Research Division\\
Faculty of Mathematics and Natural Sciences \\
Institut Teknologi Bandung \\
Jl.Ganesha 10 Bandung 40132 Indonesia}
\email{suhadi@students.itb.ac.id}
\author{E.T. Baskoro}
\address[E.T. Baskoro]{ Combinatorial Mathematics Research Division\\
Faculty of Mathematics and Natural Sciences \\
Institut Teknologi Bandung \\
Jl.Ganesha 10 Bandung 40132 Indonesia}
\address{Abdus Salam School of Mathematical Sciences\\
68-B New Muslim Town, Lahore 54600, Pakistan}
\email{ebaskoro@math.itb.ac.id}
\author{A.N.M. Salman}
\address[A.N.M. Salman]{ Combinatorial Mathematics Research Division\\
Faculty of Mathematics and Natural Sciences \\
Institut Teknologi Bandung \\
Jl.Ganesha 10 Bandung 40132 Indonesia}
\email{msalman@math.itb.ac.id}
\author{D. Suprijanto}
\address[D. Suprijanto]{ Combinatorial Mathematics Research Division\\
Faculty of Mathematics and Natural Sciences \\
Institut Teknologi Bandung \\
Jl.Ganesha 10 Bandung 40132 Indonesia}
\email{djoko@math.itb.ac.id}
\author{M. Ba\v{c}a}
\address[M. Ba\v{c}a]{Department of Appl. Mathematics, Technical University,
Letn\'{a} 9, 042 00 Ko\v{s}ice, Slovak Republic}
\email{Martin.Baca@tuke.sk}
\date{}
%\subjclass[2000]{05C12}
\dedicatory{}
\thanks{}

\begin{abstract}
A set of vertices $W$ resolves a graph $G$ if every vertex is uniquely
determined by its vector of distances to the vertices in $W$. \ A metric
dimension of $G$ is the minimum cardinality of a resolving set of $G$. A
bipartite graph $G(n,n)$ is a graph whose vertex set $V$ can be partitioned
into two subsets $V_1$ and $V_2,$ with $|V_1|=|V_2|=n,$ such that every edge
of $G$ joins $V_1$ and $V_2$. The graph $G$ is called $k$-regular if every
vertex of $G$ is adjacent to $k$ other vertices. \ In this paper, we
determine the metric dimension of $k$-regular bipartite graphs $G(n,n)$
where $k=n-1$ or $k=n-2$. \
\end{abstract}

\maketitle



\section{Introduction}

Throughout this paper, all graphs are finite, simple, and connected. \ The
\textit{vertex set} and the \textit{edge set} of graph $G$ are denoted by $%
V\left( G\right) $ and $E\left( G\right) $, respectively. \ The distance
between two distinct vertices $u,v\in V\left( G\right), $ denoted by $%
d\left( u,v\right), $ is the length of a shortest $u-v$ path in $G$. \ Let $%
W=\left\{ w_{1},\ldots ,w_{k}\right\} $ be an ordered subset of $V\left(
G\right) $. \ For $v\in V\left( G\right) $, a \textit{representation} of $v$
with respect to $W$ is defined as $k$-tuple $r\left( v \mid W\right) =\left(
d\left( v,w_{1}\right) ,\ldots ,d\left( v,w_{k}\right) \right) $. \ The set $%
W$ is a \textit{resolving set} of $G$ if every two distinct vertices $x,y\in
V\left( G\right) $ satisfy $r\left( x \mid W\right) \neq r\left( y \mid
W\right) $. \ A \textit{basis }of $G$ is a resolving set of $G$ with minimum
cardinality, and the \textit{metric dimension} of $G$ refers to its
cardinality, denoted by $\beta \left( G\right)$.

A graph $G$ is called $k$-\textit{regular} if every vertex of $G$ is
adjacent to $k$ other vertices. \ We consider a \textit{bipartite graph}
with $n$ vertices in each partitioned subset (called independent set),
denoted by $G(n,n)$.

The metric dimension in general graphs was firstly studied by Harary and
Melter \cite{Har}, and independently by Slater \cite{Sla,Slater1988}. \
Garey and Johnson \cite{Garey} showed that determining the metric dimension
of graph is NP-complete problem which is reduced from 3 dimensional matching
(3DM), while Khuller \textit{et al} \cite{Khu96} proved that it is reduced
from 3 satisfiability (3SAT). \ However, some results for certain class of
graphs have been obtained, such as cycles \cite{Char2000z}, trees \cite%
{Char2000,Har,Khu96}, fans \cite{CacePre}, wheels \cite{Buckz, CacePre,
Sha2002}, complete $n$-partite graphs \cite{Char2000,Wido2008b}, unicylic
graphs \cite{Poisson2002}, grids \cite{Melter1984}, honeycomb networks \cite%
{Manuel2008}, circulant networks \cite{Rajan}, Cayley graphs \cite{Fehr2006}%
, graphs with pendants \cite{Iswadi}, Jahangir graphs \cite{Tomescu2007},
and amalgamation of cycles \cite{IswadiPre}. \ Moreover, several researchers
have been characterized all graphs with a given metric dimension. \ Khuler
\textit{et al.} \cite{Khu96} showed that a path $P_{n}$ is the graph $G$ if
and only if $\beta \left( G\right) =1$. In \cite{Char2000}, Chartrand
\textit{et al.} showed that $G$ is $K_{n}$ if and only if $\beta \left(
G\right) =n-1$. $\ $They also proved that $\beta \left( G\right) =n-2$ if
and only if $G$ is either $K_{r,s}$ for $r,s\geq 1$, or $K_{r}+\overline{%
K_{s}}$ for $r\geq 1,s\geq 2$, or $K_{r}+\left( K_{1}\cup K_{s}\right) $ for
$r,s\geq 1$.

Study on the metric dimension of a regular graph $G\left( n,n\right) $ was
initiated by Chartrand \textit{et al}. \cite{Char2000}. They obtained the
result for $n$-regular graph $G(n,n).$ Recently, the result was generalized
to complete $k$-partite graph by Saputro \textit{et al}. \cite{Wido2008b}.
Chartrand \textit{et al} \cite{Char2000} also determined the metric
dimension of even cycle which is isomorphic to $2$-regular graph $G(n,n)$.
The purpose of this paper is to further investigate the metric dimension of
certain family of graphs, namely to determine the metric dimension of
certain regular bipartite graphs. We obtain two main results, one of them is
the following result related with an $(n-1)$-regular bipartite graph $%
G(n,n). $

\begin{theorem}
\label{Teorema1} For $n\geq 3$, if $G(n,n)$ is an $(n-1)$-regular bipartite
graph, then $\beta \left( G\right) =n-1$.
\end{theorem}

In preparing the proof for the second result we are able to obtain the
intermediate result as follows.

\begin{theorem}
\label{Teorema2} \label{Alpha_G}For $m\geq 5$, let $H$ be a connected graph
with $H=K_{m,m}\backslash E(C_{2m})$. \ Then $\beta \left( H\right)
=\left\lfloor \frac{4m}{5}\right\rfloor .$
\end{theorem}

Our second result is related with an $(n-2)$-regular bipartite graph $%
G(n,n). $ Note that every $(n-2)$-regular bipartite graph $G(n,n)$ is
isomorphic to a graph $K_{n,n}\backslash (E\left( R_{1}\right) \cup E\left(
R_{2}\right) \cup \ldots \cup E\left( R_{r}\right))$ in the theorem below.

\begin{theorem}
\label{Teorema3} For $n\geq 4$ and $r\geq 1$, let $R_{1},R_{2},\ldots ,R_{r}$
be $r$ disjoint even cycles contained in $K_{n,n}$ such that $V\left(
R_{1}\right) \cup V\left( R_{2}\right) \cup \ldots \cup V\left( R_{r}\right)
=V\left( K_{n,n}\right) $. \ For $i\in \left\{ 1,2,\ldots ,r\right\} $, let $%
G=K_{n,n}\backslash (E\left( R_{1}\right) \cup E\left( R_{2}\right) \cup
\ldots \cup E\left( R_{r}\right) )$ and $m_{i}=\frac{\left\vert V\left(
R_{i}\right) \right\vert }{2}$. \ For every $i\in \left\{ 1,2,\ldots
,r\right\} $, let $G_{i}$ be a subgraph of $G$ such that $%
G_{i}=K_{m_{i},m_{i}}\backslash E(R_{i})$. \ If $k_{1}$ is the number of
cycles $R_{i}$ where $m_{i}=2$ or $m_{i}=0$ (mod $5$), $k_{2}$ is the number
of cycles $R_{i}$ where $m_{i}=1$ (mod $5$), and $k_{3}$ is the number of
cycles $R_{i}$ where $m_{i}=2,3,$ or $4$ (mod $5$), then
\begin{equation*}
\beta (G)=%
\begin{cases}
2, & \text{if }n=4\text{,} \\
\underset{i=1}{\overset{r}{\sum }}\beta \left( G_{i}\right) , & \text{if }%
n\geq 5\text{ and }k_{1}\in \{r-1,r\}\text{ or }r=1\text{,} \\
\underset{i=1}{\overset{r}{\sum }}\beta \left( G_{i}\right) +k_{2}+k_{3}-2,
& \text{if }n\geq 5\text{, }r\geq 2\text{, }k_{1}\leq r-2\text{, and }%
k_{3}\geq 2\text{,} \\
\underset{i=1}{\overset{r}{\sum }}\beta \left( G_{i}\right) +k_{2}+k_{3}-1,
& \text{if }n\geq 5\text{, }r\geq 2\text{, }k_{1}\leq r-2\text{, and }%
k_{3}\in \{0,1\}\text{.}%
\end{cases}%
\end{equation*}
\end{theorem}

Unless otherwise stated, from now on, $G$ denotes a regular bipartite graph $%
G(n,n)$.

\section{Proof of Theorem \protect\ref{Teorema1}}

Theorem \ref{Teorema1} is a direct consequence of two lemmas in this section.

For $n\geq 3$, let $G$ be an $(n-1)$-regular bipartite graph with $%
V_{1}(G)=\left\{ x_{1},\ldots ,x_{n}\right\} $, $V_{2}(G)=\{ y_{1},\ldots
,y_{n}\} $, and $E(G)=\{x_{i}y_{j} \mid i \neq j\}$. Certainly, for $j\neq i$%
, $d(x_{i},y_{i}) =3$, $d(x_{i},y_{j}) =1$, and $d(x_{i},x_{j})=2$.

\begin{lemma}
\label{n-1L}For $n\geq 3$, let $G$ be an $(n-1)$-regular bipartite graph. \
If $W$ is a resolving set of $G$, then $W$ contains at least $n-1$ vertices.
\end{lemma}

\begin{proof}
Suppose that $W$ contains at most $n-2$ vertices. \ We define $W_{1}=$\ $%
V_{1}\cap W$ and $W_{2}=$\ $V_{2}\cap W$. Without lost of generality, let $%
|W_{1}| =p$ with $p>0$, and $|W_{2}|=q$ with $q\geq 0$. \ We define a vertex
set $A$ as the set of all vertices $x\in V_{1}\backslash W_{1}$ for which
there exist $y\in W_{2}$ such that $xy\notin E(G).$ Since $|A| \leq q$ and $%
p+q \leq n-2$, there exist two distinct vertices $x_{a},x_{b}\in
V_{1}\backslash (W_{1}\cup A) $ which are adjacent to all vertices of $W_{2}$%
. \ Then $r(x_{a} \mid W_{2}) =r(x_{b}\mid W_{2})$ which implies $%
r(x_{a}\mid W) =r(x_{b}\mid W),$ a contradiction.
\end{proof}

\begin{lemma}
\label{n-1A}For $n\geq 3$, let $G$ be an $(n-1)$-regular bipartite graph
with $V_{1}(G)=\{x_{1},\ldots ,x_{n}\} $, $V_{2}(G)=\{y_{1},\ldots ,y_{n}\} $%
, and $E(G)=\{x_{i}y_{j}\mid i \neq j\} $. \ Let $W=\{x_{1},\ldots
,x_{n-1}\} $. \ Then $W$ is a resolving set of $G$.
\end{lemma}

\begin{proof}
For distinct $i,j\in \{ 1,2,\ldots ,n-1 \},$ we have $r(y_{i}\mid W)\neq
r(y_{j}\mid W)$ since $d(x_{i},y_{i})$ $=3$ and $d( x_{i},y_{j}) =1$.
Moreover, since $r\left(y_{n}\mid W\right) =\left( 1,1,\ldots ,1\right) $
and $r\left( x_{n}\mid W\right) =\left( 2,2,\ldots ,2\right) $, we conclude
that every two different vertices $u,v \in V(G)$ satisfy $r(u\mid W)\neq
r(v\mid W)$. \ Therefore, $W$ is a resolving set of $G$.
\end{proof}

\section{Subgraph of $(n-2)$-Regular Graphs and proof of Theorem \protect\ref%
{Teorema2}}

For $n\geq 4$, let $G$ be an $(n-2)$-regular bipartite graph. Note that the
graph $G$ has a subgraph $G^{\prime }$ isomorphic to a graph obtained by
removing a hamiltonian cycle $C_{2m}$ from complete bipartite graph $K_{m,m}$
where $m\in \{2,3,...,n\}$. \ For every $v\in V(G^{\prime })\cap V_{i}(G)$
and $w\in V(G\backslash G^{\prime })\cap V_{j}(G)$ with $i,j\in \{1,2\}$, we
have $vw\in E(G),$ if $i\neq j,$ and $vw\notin E(G),$ otherwise. Therefore,
every two distinct vertices in $G^{\prime }$ must be resolved by a vertex in
$G^{\prime }.$

\begin{lemma}
\label{S_exists} For $n\geq 4$, let $G$ be an $(n-2)$-regular bipartite
graph with $\left\vert V_{1}(G)\right\vert =\left\vert V_{2}(G)\right\vert
=n $. \ Let $G^{\prime }\subseteq G$ be such that $G^{\prime
}=K_{m,m}\backslash E(C_{2m})$ with $2\leq m\leq n$. \ If $W$ is a resolving
set of $G$, then $W$ includes a basis of $G^{\prime }$.
\end{lemma}

\begin{proof}
Let $W^{\prime }$ be a basis of $G^{\prime }$. \ Suppose that there exists $%
z\in W^{\prime }$ such that $z\notin W$. \ Note that for every $v\in
V(G^\prime)\cap V_{i}(G)$ and $w\in V(G\backslash G^\prime)\cap V_{j}(G)$
with $i,j\in \{1,2\}$, we have $vw\in E(G),$ if $i \neq j,$ and $vw\notin
E(G),$ otherwise. \ Therefore, since $z\notin W$, for $i,j\in \{1,2 \}$ and $%
i\neq j$, if $z\in V_{i}$ then we obtain two possibilities below.

\begin{enumerate}
\item There exist two different vertices $x,y\in V_{j}(G)\cap V(G^{\prime })$
such that $r(x\mid W)= r(y\mid W)$.

\item There exists a vertex $x\in V_{i}\cap V(G^{\prime })$ such that $%
r(x\mid W)= r(z\mid W)$.
\end{enumerate}

In both possibilities we have a contradiction.
\end{proof}

The following lemma provides the metric dimension of a certain subgraph of $%
(n-2)$-regular bipartite graph.

\begin{lemma}
\label{2and3} For $n\geq 4$, let $G$ be an $(n-2)$-regular bipartite graph
with $\left\vert V_{1}(G)\right\vert =\left\vert V_{2}(G)\right\vert =n$. \
Let $G^{\prime }\subseteq G$ such that $G^{\prime }=K_{m,m}\backslash
E(C_{2m})$ with $m\in \{2,3,4\}$, and $W$ be a resolving set of $G$. \ \ For
$m\in \{2,3\}$, $G^{\prime }$ contributes at least $2$ vertices in $W$. \
For $m=4$, if $n=4$ then $G^{\prime }$ contributes at least $2$ vertices in $%
W$, otherwise $G^{\prime }$ contributes at least $3$ vertices in $W$.
\end{lemma}

\begin{proof}
For $m\in \{2,3,4\}$, let $V_{1}(G^{\prime })=\left\{ x_{1},\ldots
,x_{m}\right\} $, $V_{2}(G^{\prime })=\left\{ y_{1},\ldots ,y_{m}\right\} $,
and $C_{2m}=x_{1}y_{1}x_{2}y_{2}\ldots $ $x_{m}y_{m}x_{1}$. We distinguish
two cases.\newline

\textbf{Case 1:} $\beta (G^\prime)=2.$

For $n=4$ and $m=4$, we have $G^{\prime }=G$ which is isomorphic to even
cycle. \ Chartrand \textit{et al}. \cite{Char2000z} have proved that the
metric dimension of even cycle is $2$. \ Now, we assume that $n\geq 5$ and $%
m\in \{2,3\}$.

Since $\beta(G)=1$ if and only if $G$ is a path $P_n$ (\cite{Char2000}, \cite%
{Khu96}), we have $\beta (G^{\prime })\geq 2$. \ Now, we show that $\beta
\left( G^{\prime }\right) \leq 2$ by constructing a resolving set $W$ with $%
2 $ vertices. \

\begin{enumerate}
\item \textbf{For $m=2:$} We define $W=\{x_{1},y_{1}\}$. Note that, in $G$
we have $d(x_{1},y_{2})=3$ and $d(x_{1},x_{2})=2$. \ So, we obtain $%
r(x_{2}\mid W)\neq r(y_{2}\mid W)$. Therefore, $W$ is a resolving set.

\item \textbf{For $m=3:$} We define $W=\{x_{1},x_{2}\}$. \ It is easy to see
that $r(x_{3}\mid W)=(2,2)$, $r(y_{1}\mid W)=(3,3)$, $r(y_{2}\mid W)=(1,3)$,
and $r(y_{3}\mid W)=(3,1)$. \ Since every two distinct vertices $u,v$ of $%
G^{\prime }$ satisfy $r(u\mid W)\neq r(v\mid W)$, $W$ is a resolving set.
\end{enumerate}

\textbf{Case 2:} $\beta (G^\prime)=3.$

Let $m=4$ and $n \geq 5.$ First, we show that $\beta \left( G^{\prime
}\right) \leq 3$ by constructing a resolving set $W$ with $3$ vertices. \ We
define $W=\{x_{1},x_{2},x_{3}\}$. \ It is easy to see that $r(x_{4}\mid
W)=(2,2,2)$, $r(y_{1}\mid W)=(3,3,1)$, $r(y_{2}\mid W)=(1,3,3)$, $%
r(y_{3}\mid W)=(1,1,3)$, and $r(y_{4}\mid W)=(3,1,1)$. \ Since every two
distinct vertices $u,v$ of $G^{\prime }$ satisfy $r(u\mid W)\neq r(v\mid W)$%
, $W$ is a resolving set.

Next, we show that $\beta \left( G^{\prime }\right) \geq 3$. \ Suppose that $%
\beta \left( G^{\prime }\right) \leq 2$ and $S$ is the basis of $G^{\prime }$%
. \ We consider two possibilities of $S$.

\begin{enumerate}
\item $S=\{x_{i},x_{j}\}$ (or $S=\{y_{i},y_{j}\}$) where $i,j\in \{1,2,3,4\}$
and $i\neq j$

Since $m=4$, there exist two distinct vertices $x_{p},x_{q}\notin S$ (or $%
y_{p},y_{q}\notin S$). \ Since every $u\in S$ satisfies $ux_{p},ux_{q}\notin
E(G^{\prime })$, we obtain $r(x_{p}\mid S)=r(x_{q}\mid S)$, a contradiction.
\ Similarly, we obtain $r(y_{p}\mid S)=r(y_{q}\mid S)$ for $%
S=\{y_{i},y_{j}\} $.

\item $S=\{x_{i},y_{j}\}$ where $i,j\in \{1,2,3,4\}$

If $j=i$, then there exist two distinct vertices $x_{p},x_{q}\notin S$ such
that $x_{p}y_{j},x_{q}y_{j}\in E(G^{\prime })$, otherwise there exist two
distinct vertices $x_{p},x_{q}\notin S$ such that $x_{p}y_{j},x_{q}y_{j}%
\notin E(G^{\prime })$. \ Since $x_i \in S$ and $x_ix_{p},x_ix_{q}\notin
E(G^{\prime })$, both conditions imply $r(x_{p}\mid S)=r(x_{q}\mid S)$, a
contradiction.
\end{enumerate}
\end{proof}

\begin{remark}
Lemma \ref{2and3} says that the metric dimension of a subgraph $G^\prime$
above is given by
\begin{equation*}
\beta (G^{\prime })=
\begin{cases}
2, & \text{if }m\in \{2,3\} \text{ and } n \geq 5 ;\text{ or }m=4\text{ and }%
n=4, \\
3, & \text{if }m=4\text{ and }n\geq 5.%
\end{cases}%
\end{equation*}
\end{remark}

\subsection{Gap between two vertices}

Motivated by the result given in Lemma \ref{S_exists}, our observation now
is focused on a graph $H$ which is the complete bipartite graph minus its
hamiltonian cycle. Let $H=K_{m,m}\backslash E(C_{2m})$ where $m\geq 4$. \
For $m=4$, $H$ is isomorphic to an even cycle with $2m$ vertices (\cite%
{Char2000z}).

For $m\geq 5$, let $S$ be the set of two or more vertices of $H$. \ Let $%
v,w\in S$ and $P$ be a shortest $v-w$ path in $C_{2m}$. \ (It is clear that
all edges of $P$ are not elements of $E(H)$).\ We define a \emph{gap}
between $v$ and $w$ as the set of vertices in $P-\{v,w\}.$ Here, the
vertices $v$ and $w$ are called the \textit{end points}. \ The two gaps
which have at least one common end point, will be referred to as \textit{%
neighboring gaps}. \ Consequently, if $\left\vert S\right\vert =r$, then $S$
has $r$ gaps, some of gaps may be empty. \ These definitions are firstly
introduced by Buczkowski \textit{et al}. \cite{Buckz}. \ Furthermore, I.
Tomescu and I. Javaid \cite{Tomescu2007} used this gap technique to prove
the metric dimension of the Jahangir graph $J_{2n}$.

Now, let $W$ be a basis of $H$. \ We observe the following five facts.

\begin{enumerate}
\item[(i)] \emph{Every gap of $W$ contains at most four vertices.}
Otherwise, there is a gap containing at least five vertices $%
a_{1},a_{2},a_{3},a_{4},a_{5}$ of $H $ where $a_{j}a_{j+1}\notin E(H) $ with
$1\leq j\leq 4$. \ However, for every $u\in W$, $d(u,a_{2}) =d(u,a_{4}) $
which implies $r( a_{2}\mid W) =r(a_{4}\mid W)$, a contradiction.

\item[(ii)] \emph{At most one gap of $W$ contains four vertices.} Otherwise,
there exist two distinct gaps $\left\{ a_{1},a_{2},a_{3},a_{4}\right\} $ and
$\left\{ b_{1},b_{2},b_{3},b_{4}\right\} $ where $a_{j}a_{j+1},b_{j}b_{j+1}%
\notin E\left( H\right) $ for $1\leq j\leq 3$, and let $u\in W$. \ If $a_{2}$
and $b_{2}$ are in the same independent set, then $d(u,a_{2}) =d(u,b_{2})$
so $r(a_{2}\mid W) =r(b_{2}\mid W) $, otherwise $d(u,a_{2}) =d(u,b_{3})$ so $%
r(a_{2}\mid W) =r(b_{3}\mid W)$. In both cases we have a contradiction.

\item[(iii)] \emph{If a gap $A$ of $W$ contains k vertices where }$2\leq
k\leq 4$\emph{, then any neighboring gaps of $A$ contain at most one vertex.}
Otherwise, there are $k+3$ vertices $a_{1},a_{2},\ldots ,a_{k+3}$ of $H$
where $a_{j}a_{j+1}\notin E\left( H\right) $ with $1\leq j\leq k+2$, and $%
a_{k+1}$ is the only vertex among $a_{1},a_{2},\ldots ,a_{k+3}$ contained in
$W$. \ Then $r(a_{k}\mid W)=r(a_{k+2}\mid W)$, a contradiction.

\item[(iv)] \emph{For any two gaps $A$ and $B$ of $W$ containing three
vertices, both of their end points are located in different independent sets
of $H$.} Otherwise, there are ten vertices $%
a_{1},a_{2},a_{3},a_{4},a_{5},b_{1},b_{2},b_{3},b_{4},b_{5}$ of $H$ where $%
a_{j}a_{j+1},b_{j}b_{j+1}\notin E\left( H\right) $ with $1\leq j\leq 4$, and
$a_{1},a_{5},b_{1},b_{5}$ are the only vertices of $W$ from the same
independent sets of $H$. Then $r(a_{3}\mid W)=r(b_{3}\mid W)$, a
contradiction.

\item[(v)] \emph{The graph $H$ contains either gap of three vertices or gap
of four vertices.} Otherwise, there exist two distinct gaps $%
\{a_{1},a_{2},a_{3}\}$ and $\{b_{1},b_{2},b_{3},b_{4}\}$ where $%
a_{j}a_{j+1},b_{k}b_{k+1}$ $\notin E(H)$, for $1\leq j\leq 2$ and $1\leq
k\leq 3$. If $a_{2}$ and $b_{2}$ are in the same independent set of $H$,
then $r(a_{2}\mid W)=r(b_{2}\mid W)$, otherwise $r(a_{2}\mid W)=r(b_{3}\mid
W)$. We obtain a contradiction.
\end{enumerate}

Now, let $S$ be any set of vertices of $H$ which satisfies (i)-(v) above,
and let $u\in V(H)\backslash S$. There are four possibilities concerning the
position of $u$ respect to the gaps generated by $S.$

\begin{enumerate}
\item \emph{$u$ belongs to a gap of size one in $S$.} \ Let $a,b$ be
distinct end points of this gap. Then the vertex $u$ have a distance $3$ to
both $a$ and $b$, and it is the only vertex which has this distance
property. Therefore, for all $x\in V(H) $ and $x\neq u$, we have $r(x\mid S)
\neq r(u\mid S).$

\item \emph{$u$ belongs to a gap of size two in $S$.} \ Consider the
vertices $a_{1},u,a_{2},a_{3}$ of $H$ with $a_{1},a_{3}\in S$. \ Then the
vertex $u$ has a distance $3$ and $2$ to $a_{1}$ and $a_{3}$, respectively.
\ Let $v\in S\backslash \left\{ a_{1},a_{3}\right\} $. \ If $v$ is in the
same independent set as $u$ then $d\left( u,v\right) =2$, otherwise $d\left(
u,v\right) =1$. \ By observation (iii), the vertex $u$ is the only one
having all of these distance properties. \ Therefore, for all $x\in V\left(
H\right) $ and $x\neq u$, we have $r(x\mid S)\neq r(u\mid S)$.

\item \emph{$u$ belongs to a gap of size three in $S$.} \ Consider the
vertices $a_{1},a_{2},a_{3},a_{4},a_{5}$ of $H$ with only $a_{1},a_{5}\in S$%
. \ If $u=a_{2}$, then $d\left( u,a_{1}\right) =3$, and if $u=a_{3}$, then $%
d\left( u,a_{1}\right) =2$. \ Let $v\in S\backslash \left\{
a_{1},a_{5}\right\} $. \ If $v$ is in the same independent set as $u$ then $%
d(u,v)=2$, otherwise $d(u,v)=1$. \ By observations (iii)-(v), the vertex $u$
is the only one having all of these distance properties. \ Therefore, for
all $x\in V(H)$ and $x\neq u$, we have $r(x\mid S)\neq r(u\mid S)$.

\item \emph{$u$ belongs to a gap of size four in $S$.} \ Consider the
vertices $a_{1},a_{2},a_{3},a_{4},a_{5},a_{6}$ of $H$ with only $%
a_{1},a_{6}\in S$. \ If $u=a_{2}$, then $d(u,a_{1})=3$ and $d\left(
u,a_{6}\right) =2$. \ If $u=a_{3}$, then $d(u,a_{1})=2$ and $d(u,a_{6})=1$.
\ Let $v\in S\backslash \{a_{1},a_{6}\}$. If $v$ is in the same independent
set as $u$ then $d(u,v)=2$, otherwise $d(u,v)=1$. By observations (i)-(iii)
and (v), the vertex $u$ is the only one having all of these distance
properties. Therefore, for all $x\in V(H)$ and $x\neq u$, we have $r(x\mid
S)\neq r(u\mid S)$.
\end{enumerate}

Consequently, any set $S$ having properties (i)-(v) resolves $V(H)$. \ Now,
we provide a proof of Theorem \ref{Teorema2}.

\subsection{Proof of Theorem \protect\ref{Teorema2}}

For $m\geq 5$, let $V_{1}(H)=\left\{ x_{1},\ldots ,x_{m}\right\} $, $%
V_{2}(H)=\left\{ y_{1},\ldots ,y_{m}\right\} $, and $%
C_{2m}=x_{1}y_{1}x_{2}y_{2}\ldots $ $x_{m}y_{m}x_{1}$.

\textbf{Claim 1:} $\beta \left( H\right) \leq \left\lfloor \frac{4m}{5}%
\right\rfloor.$

We show that $\beta \left( H\right) \leq \left\lfloor \frac{4m}{5}%
\right\rfloor $ by constructing a resolving set $W$ with $\left\lfloor \frac{%
4m}{5}\right\rfloor $ vertices. \ We consider the integer $k\geq 1$. \ We
obtain four cases as follows.

\begin{enumerate}
\item[(a)] $m=0$ or $1$ (mod $5$)\

Let $m=5k$ or $m=5k+1$ with $k\geq 1$ (hence, $\left\lfloor \frac{4m}{5}%
\right\rfloor =4k$). \ We define $W=\{y_{5j+1},y_{5j+2},$ $%
x_{5j+4},x_{5j+5}\mid 0\leq j\leq k-1\}$. \ Since $W$ contains $4k$ vertices
and satisfies (i)-(v), then $W$ is a resolving set.

\item[(b)] $m=2$ (mod $5$)

Let $m=5k+2$ with $k\geq 1$ (hence, $\left\lfloor \frac{4m}{5}\right\rfloor
=4k+1$). \ We define $W=\{y_{5j+1},y_{5j+2},$ $x_{5j+4},x_{5j+5}\mid 0\leq
j\leq k-1\}\cup \{x_{5k+1}\}$. \ Since $W$ contains $4k+1$ vertices and
satisfies (i)-(v), then $W$ is a resolving set.

\item[(c)] $m=3$ (mod $5$)

Let $m=5k+3$ with $k\geq 1$ (hence $\left\lfloor \frac{4m}{5}\right\rfloor
=4k+2$). \ We define $W=\{y_{5j+1},y_{5j+2},$ $x_{5j+4},$ $x_{5j+5}\mid
0\leq j\leq k-1\}\cup \{y_{5k+1},y_{5k+2}\}$. \ Since $W$ contains $4k+2$
vertices and satisfies (i)-(v), then $W$ is a resolving set.

\item[(d)] $m=4$ (mod $5$)

Let $m=5k+4$ with $k\geq 1$ (hence $\left\lfloor \frac{4m}{5}\right\rfloor
=4k+3$). \ We define $W=\{y_{5j+1},y_{5j+2},$ $x_{5j+4},x_{5j+5}\mid 0\leq
j\leq k-1\}\cup \{y_{5k+1},$ $y_{5k+2},y_{5k+3}\}$. \ Since $W$ contains $%
4k+3$ vertices and satisfies (i)-(v), then $W$ is a resolving set.
\end{enumerate}

\textbf{Claim 2:} $\beta \left( H\right) \geq \left\lfloor \frac{4m}{5}%
\right\rfloor.$

Let $S$ be a basis of $H$. \ We consider two cases as follows.

\begin{enumerate}
\item $\beta \left( H\right) $ is even.

Let $\left\vert S\right\vert =2l$ for some integer $l\geq 2$. \ By (iii),
there are at most $l$ gaps which contain more than one vertex. \ By (i),
(ii), (iv), and (v), all of them contain $2$ vertices, except possibly
either two contain $3$ vertices or one contains $4$ vertices. \ Then, the
number of vertices belonging to the gaps of $S$ is at most $3l+2$. \ Hence $%
2m-2l\leq 3l+2$, which implies $\left\vert S\right\vert =2l\geq \left\lceil
\frac{4m}{5}-\frac{4}{5}\right\rceil \geq \left\lfloor \frac{4m}{5}%
\right\rfloor $.

\item $\beta \left( H\right)$ is odd.

Let $\left\vert S\right\vert =2l+1$ for some integer $l\geq 2$. \ By (iii),
there are at most $l$ gaps which contain more than one vertex. \ By (i),
(ii), (iv), and (v), all of them contain $2$ vertices, except possibly
either two contain $3$ vertices or one contains $4$ vertices. \ Then, the
number of vertices belonging to the gaps of $S$ is at most $3l+3$. \ Hence $%
2m-2l\leq 3l+3$, which implies $\left\vert S\right\vert =2l+1\geq
\left\lceil \frac{4m}{5}-\frac{6}{5}+1\right\rceil \geq \left\lfloor \frac{4m%
}{5}\right\rfloor $.
\end{enumerate}

\qed

\section{Gap in $(n-2)$-regular bipartite graph and proof of Theorem \protect
\ref{Teorema3}}

For $n\geq 4$, we consider certain cycles contained in a complete bipartite
graph $K_{n,n}.$ For $r\geq 1$, let $R_{1},R_{2},\ldots ,R_{r}$ be $r$
disjoint even cycles contained in $K_{n,n}$ such that $V\left( R_{1}\right)
\cup V\left( R_{2}\right) \cup \ldots \cup V\left( R_{r}\right) =V\left(
K_{n,n}\right) $. \ Then $K_{n,n}\backslash (E\left( R_{1}\right) \cup
E\left( R_{2}\right) \cup \ldots \cup E\left( R_{r}\right) )$ is an $(n-2)$%
-regular bipartite graph $G(n,n)$.

Let $H=K_{n,n}\backslash (E\left( R_{1}\right) \cup E\left( R_{2}\right)
\cup \ldots \cup E\left( R_{r}\right) )$ and $m_{i}=\frac{\left\vert V\left(
R_{i}\right) \right\vert }{2}$. \ For every $i\in \left\{ 1,2,\ldots
,r\right\} $, let $G_{i}$ be a subgraph of $H$ such that $%
G_{i}=K_{m_{i},m_{i}}\backslash E(R_{i})$. \ So, $G_{i}=K_{m_{i},m_{i}}%
\backslash E(C_{2m_{i}})$. \ \ \

\begin{figure}[h]
\includegraphics[bb=0 0 230 280, scale=0.55]{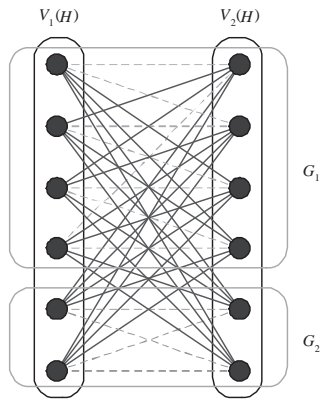}
\caption{An example of graph $H$}
\label{fig1}
\end{figure}

Now, we consider a basis $W$ of $H$. \ By Lemma \ref{S_exists}, $W$ is a
union of the resolving sets of $G_{i}$ with $1\leq i\leq r$. \ So, we have $%
\beta (H)\geq \beta (G_{1})+\beta (G_{2})+...+\beta (G_{r})$. \ However,
there are three conditions that must be satisfied by $W$ for $r\geq 2$,
deduced from (ii), (iv), and (v), respectively.

\begin{enumerate}
\item[(a)] At most one gap in $W$ contains four vertices. \

\item[(b)] For any two gaps $A$ and $B$ of $W$ containing three vertices,
both of their end points are located in different independent sets of $H$. \

\item[(c)] $H$ contains either gap of three vertices or gap of four
vertices. \
\end{enumerate}

Although the basis of $G_{i}$ with $1\leq i\leq r$ are resolving sets, there
may exists some $i\in \left\{ 1,2,\ldots ,r\right\} $ such that a resolving
set of $G_{i}$ which is contained in $W$, is not the basis of $G_{i}$. \ For
example, if there exists $j\in \left\{ 1,2,\ldots ,r\right\} $ and $j\neq i$
such that the basis of both $G_{i}$ and $G_{j}$ have a gap containing four
vertices, then by (a), we cannot use the basis of $G_{j}$ in $G$. \ We must
add at least one more vertex on the basis of $G_{j}$ such that the new
resolving set of $G_{j}$ satisfies (a)-(c). \ So, we need to know the gaps
property of the basis of $G_{i}=K_{m_{i},m_{i}}\backslash E(C_{2m_{i}})$,
which can be seen in Lemmas \ref{1234mod5}, \ref{0mod5}, \ref{1mod5}, and %
\ref{234mod5}. \

\begin{lemma}
\label{1234mod5}For $n\geq 5$, let $G(n,n)$ be an $(n-2)$-regular bipartite
graph. Let $G^{\prime }\subseteq G$ such that $G^{\prime }=K_{m,m}\backslash
E(C_{2m})$ with $m\in \{3,4,...,n\}$. \ If $m=1,2,3,$ or $4$ $($mod $5)$,
then the basis of $G^{\prime }$ has a gap containing at least three vertices.
\end{lemma}

\begin{proof}
Suppose that $G^{\prime }$ only has the gaps which contain at most two
vertices. \ We distinguish four cases:

\noindent \textbf{Case 1:} $m=1$ (mod $5$).

Let $m=5k+1$ with $k\geq 1$. \ By (iii), there are at most $\frac{\beta
\left( G^{\prime }\right) }{2}$ gaps which contain two vertices each and at
most $\frac{\beta \left( G^{\prime }\right) }{2}$ gaps which contain one
vertex each. \ Then, $\left\vert V\left( G^{\prime }\right) \right\vert \leq
\frac{5}{2}\beta \left( G^{\prime }\right) $. \ By Theorem \ref{Alpha_G}, we
have $\left\vert V\left( G^{\prime }\right) \right\vert \leq 10k$. \ Since
there are $10k+2$ vertices, we obtain a contradiction.

\noindent \textbf{Case 2:} $m=2$ (mod $5$).

Let $m=5k+2$ with $k\geq 1$. \ By (iii), there are at most $\frac{\beta
\left( G^{\prime }\right) -1}{2}$ gaps which contain two vertices each and $%
\frac{\beta \left( G^{\prime }\right) +1}{2}$ gaps which contain one vertex
each. \ Then, $\left\vert V\left( G^{\prime }\right) \right\vert \leq \frac{5%
}{2}\beta \left( G^{\prime }\right) -\frac{1}{2}$. \ By Theorem \ref{Alpha_G}%
, we have $\left\vert V\left( G^{\prime }\right) \right\vert \leq 10k+2$. \
Since there are $10k+4$ vertices, we obtain a contradiction.

\noindent \textbf{Case 3:} $m=3$ (mod $5$).

Let $m=5k+3$ with $k\geq 0$. \ By (iii), there are at most $\frac{\beta
\left( G^{\prime }\right) }{2}$ gaps which contain two vertices each and $%
\frac{\beta \left( G^{\prime }\right) }{2}$ gaps which contain one vertex
each. \ Then, $\left\vert V\left( G^{\prime }\right) \right\vert \leq \frac{5%
}{2}\beta \left( G^{\prime }\right) $. \ By Lemma \ref{2and3} for $m\in
\{3,4\}$ or Theorem \ref{Alpha_G} otherwise, we have $\left\vert V\left(
G^{\prime }\right) \right\vert \leq 10k+5$. \ Since there are $10k+6$
vertices, we obtain a contradiction.

\noindent \textbf{Case 4:} $m=4$ (mod $5$).

Let $m=5k+4$ with $k\geq 0$. \ By (iii), there are at most $\frac{\beta
\left( G^{\prime }\right) -1}{2}$ gaps which contain two vertices each and $%
\frac{\beta \left( G^{\prime }\right) +1}{2}$ gaps which contain one vertex
each. \ Then, $\left\vert V\left( G^{\prime }\right) \right\vert \leq \frac{5%
}{2}\beta \left( G^{\prime }\right) -\frac{1}{2}$. \ By Theorem \ref{Alpha_G}%
, we have $\left\vert V\left( G^{\prime }\right) \right\vert \leq 10k+7$. \
Since there are $10k+8$ vertices, we obtain a contradiction.
\end{proof}

\begin{lemma}
\label{0mod5} For $n\geq 4$, let $G(n,n)$ be an $(n-2)$-regular bipartite
graph. \ Let $G^{\prime }\subseteq G$ such that $G^{\prime
}=K_{m,m}\backslash E(C_{2m})$ with $m\in \{2,3,...,n\}$. \ If $m=2$ or $m=0$
$($mod $5)$, then there exists a basis of $G^{\prime }$ where every gap
contains at most two vertices.
\end{lemma}

\begin{proof}
For $m\geq 2$, let $V_{1}(H)=\left\{ x_{1},\ldots ,x_{m}\right\} $, $%
V_{2}(H)=\left\{ y_{1},\ldots ,y_{m}\right\} $, and $%
C_{2m}=x_{1}y_{1}x_{2}y_{2}\ldots $ $x_{m}y_{m}x_{1}$. \ We distinguish two
cases:

\noindent \textbf{Case 1:} $m=2.$

Note that, if $n=4$, then it is impossible to have $G^{\prime }$ with $m=2$
since $G$ is a disconnected graph. \ Now, we assume that $n\geq 5$. \ We
define $S=\left\{ x_{1},y_{1}\right\} $. \ By (iii), we obtain that $S$ is a
resolving set of $G^{\prime }$. \ Since $\left\vert S\right\vert =\beta
\left( G^{\prime }\right) $, $S$ is a basis of $G^{\prime }$.

\noindent \textbf{Case 2:} $m=0$ (mod $5$).

For $m\geq 5$ and the integer $k\geq 1$, let $m=5k$. \ We define $%
S=\{y_{5j+1},y_{5j+2},$ $x_{5j+4},x_{5j+5}\mid0\leq j\leq k-1\}$. \ It is
easy to see that every gap of $S$ contains at most two vertices. \ By (iii),
we obtain that $S$ is a resolving set of $G^{\prime }$. \ Since $\left\vert
S\right\vert =\beta \left( G^{\prime }\right) $, $S$ is a basis of $%
G^{\prime }$.
\end{proof}

Lemma \ref{1234mod5} says that the basis of subgraph $G^\prime$ of $G(n,n)$
has a gap containing at least three vertices, for $m=1,2,3,4 ~(\text{mod }%
5). $ The following two lemmas are a kind of special case of Lemma 5. Lemma %
\ref{1mod5} shows that a basis of $G^\prime$ has at least two gaps
containing three vertices, for $m=1~(\text{mod }5),$ while for $m=2,3,4~(%
\text{mod }5),$ Lemma \ref{234mod5} shows the existence of a basis of $%
G^\prime$ which has one gap containing exactly three vertices.

\begin{lemma}
\label{1mod5}For $n\geq 4$, let $G(n,n)$ be an $(n-2)$-regular bipartite
graph. \ Let $G^{\prime }\subseteq G$ such that $G^{\prime
}=K_{m,m}\backslash E(C_{2m})$ with $m\in \{3,4,...,n\}$. \ If $m=1$ $($mod $%
5)$, then a basis of $G^{\prime }$ has at least $2$ gaps containing three
vertices.
\end{lemma}

\begin{proof}
Suppose that $G^{\prime }$ has at most one gap containing three vertices. \ $%
G^{\prime }$ cannot have a gap containing four vertices, by (v). \ So, the
other gaps contain at most two vertices. \ Let $m=5k+1$ with $k\geq 1$. \ By
(iii), there are at most $\frac{\beta \left( G^{\prime }\right) }{2}-1$ gaps
which contain two vertices each and $\frac{\beta \left( G^{\prime }\right) }{%
2}$ gaps which contain one vertex each. \ Then, $\left\vert V\left(
G^{\prime }\right) \right\vert \leq \frac{5}{2}\beta \left( G^{\prime
}\right) +1$. \ By Theorem \ref{Alpha_G}, we have $\left\vert V\left(
G^{\prime }\right) \right\vert \leq 10k+1$. \ Since there are $10k+2$
vertices, we obtain a contradiction.
\end{proof}

\begin{lemma}
\label{234mod5}For $n\geq 4$, let $G(n,n)$ be an $(n-2)$-regular bipartite
graph. \ Let $G^{\prime }\subseteq G$ such that $G^{\prime
}=K_{m,m}\backslash E(C_{2m})$ with $m\in \{3,4,...,n\}$. \ If $m=2,3,$ or $%
4 $ $($mod $5)$, then there exists a basis of $G^{\prime }$ which contains
one gap of three vertices.
\end{lemma}

\begin{proof}
Let $V_{1}(H)=\left\{ x_{1},\ldots ,x_{m}\right\} $, $V_{2}(H)=\left\{
y_{1},\ldots ,y_{m}\right\} $, and $C_{2m}=x_{1}y_{1}x_{2}y_{2}\ldots $ $%
x_{m}y_{m}x_{1}$.

For $m=3$, we define $W=\left\{ y_{1},y_{2}\right\} $. \ Since $|W|=\beta
\left( G^{\prime }\right) $ and satisfies (i)-(v), then $W$ is a basis of $%
G^{\prime }$.

For $m=4$, we define $W=\left\{ y_{1},y_{2},y_{3}\right\} $. \ Since $%
|W|=\beta \left( G^{\prime }\right) $ and satisfies (i)-(v), then $W$ is a
basis of $G^{\prime }$.

For $m\geq 5$ and the integer $k\geq 1$, we consider the following $3$ cases.

\begin{enumerate}
\item $m=2$ (mod $5$)

Let $m=5k+2$ with $k\geq 1$. \ We define $W=%
\{y_{5j+1},y_{5j+2},x_{5j+4},x_{5j+5}$ $\mid0\leq j\leq k-1\}\cup
\{x_{5k+1}\}$. \

\item $m=3$ (mod $5$)

Let $m=5k+3$ with $k\geq 1$. \ We define $W=%
\{y_{5j+1},y_{5j+2},x_{5j+4},x_{5j+5}$ $\mid0\leq j\leq k-1\}\cup
\{y_{5k+1},y_{5k+2}\}$. \

\item $m=4$ (mod $5$)

Let $m=5k+4$ with $k\geq 1$. \ We define $W=%
\{y_{5j+1},y_{5j+2},x_{5j+4},x_{5j+5}$ $\mid0\leq j\leq k-1\}\cup
\{y_{5k+1},y_{5k+2},y_{5k+3}\}$. \
\end{enumerate}

It is easy to see that $W$ from all three cases above contains one gap of
three vertices. \ Since $|W| =\beta \left( G^{\prime }\right) $ and
satisfies (i)-(v), then $W$ is a basis of $G^{\prime }$.
\end{proof}

We are now ready to prove Theorem \ref{Teorema3}.

\subsection{Proof of Theorem \protect\ref{Teorema3}}

For $n=4$, there is only one possibility for $G,$ namely $G$ is isomorphic
with the even cycle graph. \ Chartrand \textit{et al.} \cite{Char2000z}
proved that the metric dimension of the even cycle graph is equal to $2$. \
Now, we assume that $n\geq 5$.

The second case for $\beta (G)$ is a direct consequence of Theorem \ref%
{Teorema2}, Lemma \ref{0mod5}, and conditions (a)-(c) above.

For the last two cases, let $G^{\prime }=K_{m_{i},m_{i}}\backslash E(R_{i})$
be a subgraph of $G$ and $W^{\prime }$ be a basis of $G^{\prime }$. Suppose $%
G^{\prime}$ has a gap containing three vertices $a_{1},a_{2},a_{3}$ where $%
a_{j}a_{j+1}\notin E\left( G^{\prime }\right) $ with $1\leq j\leq 2$, or
containing four vertices $a_{1},a_{2},a_{3},a_{4}$ of $G^{\prime }$ where $%
a_{j}a_{j+1}\notin E\left( G^{\prime }\right) $ with $1\leq j\leq 3$. \ It
is easy to see that $W^{\ast }=W^{\prime }\cup \{a_{2}\}$ is a resolving set
of $G^{\prime }$ which all the gaps contain at most two vertices. \ So, by
using this property, Theorem \ref{Alpha_G}, Lemma \ref{S_exists} - Lemma \ref%
{234mod5}, and also (a)-(c), we prove the last two cases. \qed \newline

\noindent
{\bf Acknowledgement}.

The authors are thankful to the anonymous referee for some comments that
helped to improve the presentation of the manuscript.

\end{document}